\documentstyle[12pt]{article}

\newcommand{\sect}[1]{\setcounter{equation}{0}\section{#1}}

\newcommand{\f}{\frac}

\newcommand{\p}{\partial}
\newcommand{\FR}{Fun(\rn_q^N)}
\newcommand{\DFR}{Dif\!f(\rn_q^N)}
\newcommand{\FS}{Fun(SO_q(N))}

\newcommand{\bp}{\bar{\partial}}

\newcommand{\ot}{\otimes}
\newcommand{\La}{\Lambda}

\newcommand{\bc}{\begin{center}}
\newcommand{\ec}{\end{center}}
\newcommand{\be}{\begin{equation}}
\newcommand{\ee}{\end{equation}}

\def\lcross{{>\!\!\!\triangleleft}}

\font\bb=msym10
\newcommand{\cn}{\mbox{\bb C}}
\newcommand{\rn}{\mbox{\bb R}}

\newtheorem{prop}{Proposition}

\normalbaselines
\begin{document}
\begin{center}
\vspace*{1.0cm}

{\LARGE{\bf $SO_q(N)$-isotropic Harmonic Oscillator on the Quantum Euclidean
Space $\rn_q^N$}
\footnote{To appear in the proceedings of the Clausthal Symposium on
 ``Nonlinear, Dissipative, Irreversible Quantum Systems'', August 1994.}}

\vskip 1.5cm

{\large {\bf Gaetano Fiore\footnote{Alexander-von-Humboldt fellow}}}

\vskip 0.5 cm

Ls. Prof. Wess, Sektion Physik Universit\"at M\"unchen \\
Theresienstrasse 37, D-80333 M\"unchen \& Germany \\
LMU-TPW 94-26, q-alg/9501023, December 1994
\end{center}

\vspace{0.2 cm}

\begin{abstract}
We briefly describe the construction of a consistent $q$-deformation of the
quantum mechanical
isotropic harmonic oscillator on ordinary $\rn^N$ space.
\end{abstract}

\vspace{1 cm}

\section*{Introduction}

{}~~~In this report we briefly describe how
to formulate a consistent quantum-mechanical one-particle system
on a ``noncommutative manifold'' (the quantum Euclidean space $\rn_q^N$) with a
 ``non-classical''
space symmetry (the one carried by the quantum group $SO_q(N)$). This system is
a $q$-deformation of the isotropic harmonic oscillator on ordinary
$\rn^N$-space; $q\rightarrow 1$ is the corresponding ``classical limit''.
The guiding idea is to mimic  in a $q$-deformed setting
the ordinary realization of one-particle quantum mechanics in
configuration space. Essential references are \cite{car,fio1,fio2,fio4,fioth},
to which we refer the reader for details and proofs.

Incidentally, the problem considered here is rather different from that of the
q-deformed harmonic oscillators treated by other authors \cite{bie}: there,
creation/annihilation operators
with some prescribed commutations relation are postulated from the very
beginning without any reference to a geometrical framework. The deformation
considered there concerned the well-known hidden $su(n)$ symmetries of the
harmonic oscillator hamitonians. Here,
on the contrary,
the deformation concerns the rotation symmetry of the space itself
and creation/destruction operators are constructed out of the
deformed `` coordinates '' and `` derivatives ''.

The general context for a proper location of the work is a research program
\cite{wess1,majgr} consisting of:
1) a continuous deformation
of both space(time) geometry and its related fundamental symmetries through
socalled
quantum (and/or braided) spaces/groups (which are characteristic examples of
noncommutative
geometries); 2) the formulation of quantum mechanics and quantum field theory
on them.
The main physical motivation for such a program is that it represents a radical
approach to
long-standing problems in QFT, e.g. ultraviolet divergencies, quantizing
gravity, etc, since one
modifies the microscopic structure of spacetime. Of
course the effects of such a modification should be
practically undetectable at least in the domain
of observability of already well-explained physical phenomena.

The approach is an algebraic one which seems in deep agreement with the spirit
of quantum mechanics.
In fact, one of the essential features of the latter is
that it invalidates the classical geometrical
description of the state of
a physical system as a point in the corresponding phase space; the naive
formulation of this fact are Heisenberg's indetermination relations.
The notion of `` points '' in continuous phase space is replaced
by the more general notion of a
(${\bf C}^*$-)algebra of operators acting on a Hilbert space; consequently, the
geometry of the phase space is broken.
Geometry may be recovered as a useful structure for the spectral
decomposition of operators. For instance, the spectral decomposition of the
vector components of the position operator of one quantum particle
typically is realized on a manifold coinciding with classical
configuration space.
This is due to the axiom that these components commute, $[x^i,x^j]=0$.
In a noncommutative-geometric approach to quantum mechanics one essentially
releases the latter axiom. One would be interested in a tiny deformation of
these commutation
relations, but of course this can be done in infinitely too many ways.

Requiring that the deformed relations keep track of the spacetime symmetries
puts severe
constraints to the possible deformations. Quantum spaces and related quantum
groups have been conceived to satisfy
them, in that they are {\it coupled} deformations of space(time) geometries
and their fundamental symmetries.

\sect{Preliminaries}

{}~~~$\hat R_q:=||\hat R^{ij}_{hk}||$ will denote the braid matrix
of the quantum group $SO_q(N)$ \cite{frt}, $C:=||C_{ij}||$ the corresponding
q-deformed metric matrix. Here $i,j,h,j$ belong to
$\{-n,-n+1,...,-1,0,1,...n\}$ if
$N=2n+1$, and to $\{-n,-n+1,...,-1,1,...n\}$ if $N=2n$.
Indices are raised
and lowered through the metric matrix $C$, e.g.
$a_i=C_{ij}a^j,$ $a^i=C^{ij}a_j$;
$C$ is not symmetric an coincides with its inverse, $C^{ij}=C_{ij}$.

Both $C$ and $\hat R$
depend on $q$ and are real for
$q\in {\bf R}$; explicit expressions can be found in Ref. \cite{frt}.
$\hat R$ admits the very useful decomposition
\be
\hat R_q = q {\cal P}_S - q^{-1} {\cal P}_A +q^{1-N}{\cal P}_1.
\ee
${\cal P}_S,{\cal P}_A,{\cal P}_1$ are the projection operators
onto the three eigenspaces of $\hat R$ (the latter have respectively dimensions
$\f{{N(N+1)}}{ 2}-1,\f{{N(N-1)}}{ 2},1$): they project the tensor
product $x\otimes x$ of the fundamental corepresentation $x$ of $SO_q(N)$
into the corresponding irreducible corepresentations (the symmetric modulo
trace, antisymmetric and trace, namely the q-deformed
versions of the corresponding ones of $SO(N)$). The projector ${\cal P}_1$ is
related to the metric matrix $C$ by ${\cal P}_{1~hk}^{~~ij}=\f{C^{ij}C_{hk}}{
Q_N}$; the factor $Q_N$ is defined by $Q_N:=C^{ij}C_{ij}$, and
$Q_N\rightarrow N$.
$\hat R^{\pm 1},C$ satisfy the relations
\be
[f(\hat R), P\cdot (C\otimes C)]=0~~~~~~~~~
f(\hat R_{12})\hat R^{\pm 1}_{23}\hat R^{\pm 1}_{12}=
\hat R^{\pm 1}_{23}\hat R^{\pm 1}_{12}f(\hat R_{23})
\ee
($P$ is the permutator: $P^{ij}_{hk}:=\delta^i_k\delta^j_h$ and $f$ is any
polynomial function);
in particular this holds for $f(\hat R)=\hat R^{\pm 1},{\cal P}_A,{\cal P}_S,
{\cal P}_1$.

{}~~~The unital
algebra $Dif\!f({\bf R}_q^N)$ of differential operators on the real quantum
euclidean space ${\bf R}_q^N$ is defined as the space of formal series in the
(ordered)
powers of the $\{x^i\},\{\p_i\}$ variables (the $q$-deformed coordinates and
derivatives,
respectively) with complex coefficients, modulo the commutation relations
\be
{\cal P}_{A~hk}^{~~ij}x^h x^k =0.     ~~~~~~~~~~~~~
{\cal P}_{A~hk}^{~~ij}\p^h \p^k =0.
\ee
and the derivation relations
\be
\p_i x^j = \delta^j_i+q\hat R^{jh}_{ik} x^k\p_h.
\ee
The subalgebra $Fun({\bf R}_q^N)$ of `` functions '' on ${\bf R}_q^N$ is
generated by $\{x^i\}$ only. We only give as an example relations (1.3) in the
case N=3.
They amount to the three independent relations
\be
x^{-1}x^0-qx^0x^{-1}=0~~~~~~~x^0x^1-qx^1x^0=0~~~~~~~x^1x^{-1}-x^{-1}x^1=
(q^{\f 12}-q^{-\f 12})x^0x^0;
\ee
in the limit $q\rightarrow 1$ we get back commuting coordinates and therefore
ordinary geometry.
More details can be found in Ref. \cite{car} and  in the contribution of
U. Carow-Watamura $\&$ S. Watamura to these
proceedings.

{}~~~If $q\in {\bf R}^+$ one can introduce an antilinear involutive
antihomomorphism
$*$:
\be
*^2=id~~~~~~~~~~~~~~~~~~~~~~~(AB)^*=B^*A^*
\ee
on $Dif\!f({\bf R}_q^N)$ (the ``complex conjugation'').  On the basic variables
$x^i$  $*$ is
defined by
\be
(x^i)^*=x^jC_{ji}
\ee
whereas the complex conjugates of the derivatives $\p^i$ are not
combinations of the derivatives themselves. It is useful to introduce
barred derivatives $\bar \p^i$ through
\be
(\p^i)^*=-q^{-N}\bar \p^j C_{ji}.
\ee
They satisfy relation $(1.3)_2$ and the analog of (1.4) with
$q,\hat R$ replaced by $q^{-1},\hat R^{-1}$. These $\bar\p$ derivatives can
be expressed as nonlinear functions of $x,\p$ \cite{og}.

{}~~~By definition a $\FS$-scalar $I(x,\p)\in Dif\!f({\bf R}_q^N)$ transforms
trivially
under the coaction associated to the quantum group of symmetry
$SO_q(N,{\bf R})$ \cite{frt}. Any q-scalar polynomial
$I(x,\p)\in Dif\!f({\bf R}_q^N)$
of degree $2p$ in $x,\p$
can be expressed as an ordered polynomial in two particular q-scalar
variables (see for instance Appendix C of [12]),
namely the square lenght  $xCx:= x^iC_{ij}x^j$ and the laplacian
$\Delta:= \p^iC_{ij}\p^j$.
In Ref. \cite{fio4} the Hopf algebra $U_q(so(N))$ \cite{frt} was realized as
the subalgebra
of $\DFR$ characterized by the condition that its elements commute with the
$\FS$-scalar elements of $\DFR$.
One set of generators of such a subalgebra is $\{l^{ij}\}_{i\neq j}$, where
\be
l^{ij}:={\cal P}_{A~hk}^{~~ij}x^h\p^k\La^{-1}
\ee
and $\La$ is the ``dilaton'' defined by $\La x^i=q x^i\La$,
$\La\p^i=q^{-1}\p^i\La$.
Each $l^{ij}$ can be interpreted as a $q$-deformed angular momentum component,
since it commutes with $q$-scalars and in the classical limit reduces to the
ordinary angular momentum component generating a rotation in the $(i,j)$ plane
of $\rn^N$, $l^{ij}\stackrel{q\rightarrow 1}{\longrightarrow}
x^i\p^j-x^j\p^i$. $l\cdot l:=l^{ij}l_{ji}$ itself is a scalar and  therefore
commutes with each $l^{hk}$; it is a quadratic function of the quadratic
casimir of $U_q(so(N))$
and will be called the ``square angular momentum''. A linear function $B$ of
the quadratic casimir is
\be
B:=\La^{-1}(1+\f{q^2-1}{1+q^{2-N}}x^i\p_i);
{}~~~~~~\Rightarrow~~~~~~~
1=(B)^2-\f{(q^2-1)(q^2-q^{-2})}{(1+q^{2-N})(1+q^{N-4})}
(l\cdot l).
\ee

Incidentally, one can extend this realization of $U_q(so(N))$ by adding
$q$-derivatives as
generators of $q$-translations, in such a way to represent a Euclidean Hopf
algebra $U_q(e^N)$ \cite{fioeu}.

{}~~~Riemann integration $\int dV$ over the Euclidean space $\rn^N$ is
covariant
under the
action of the Euclidean group $E^N:=\rn^N\lcross SO(N)$, and in particular is
invariant under
finite translations. In infinitesimal form, the latter invariance implies
the validity of Stoke's theorem for all integrable functions
which are differentiable. $E^N$-covariance of the Riemann integration is
essential
in allowing a $E^N$-covariant description of ordinary physics. Similarly,
in view of the formulation of $E^N_q$-covariant Physics,
in Ref. \cite{fio2,fioth}
a $q$-deformed integration $\int d_qV$ on $\rn_q^N$ was constructed by
requiring its covariance
w.r.t. the quantum Euclidean group $E_q^N:=\rn_q^N\lcross SO_q(N)$
\cite{maj2,schl}. Actually,
its invariance w.r.t. to finite $q$-translations (or equivalently
the validity of $q$-Stoke's theorem, in infinitesimal form), together
with the obvious requirement of linearity, are enough to allow its
{\it construction}.
 Its $SO_q(N,\bf R)$-covariance
follows from the $SO_q(N,\bf R)$-covariance of the
``braided coaddition'' \cite{maj2} (or, equivalently, of the differential
calculus \cite{car} on $\rn_q^N$).

In the classical case, if $f=P_n(x)exp[-a|x|^2]$
($P_n$ denotes a polynomial of degree $n$ in $x$ and $|x|^2$ the square
lenght), then
\be
\p^i P_n(x)exp[-a|x|^2]|=  P_{n-1}(x)exp[-a|x|^2]+
 P_{n+1}exp[-a|x|^2];
\ee
Stoke's theorem then implies
\be
\int dV~P_{n-1}(x)exp[-a|x|^2]+ \int dV~P_{n+1}exp[-a|x|^2]=0.
\ee
This relation allows to recursively define
the integral $\int dV~f$ (for any function $f$ of the same kind)
in terms of $\int dV~exp[-a|x|^2]$
(which fixes the normalization of the integration);
moreover, since the space of functions $\{P(x)exp[-a|x|^2]\}$ is dense in
${\cal L}^1(\rn^N)$,
${\cal L}^2(\rn^N)$, one can approximate as much as desired integrals of other
functions in this
way. This situation can be summarized by saying that
$exp[-a|x|^2]$ can be taken as a ``reference function'' for the algebraic
construction
of Riemann integration.
In an analogous way, $\int d_qV$ was algebraically {\it constructed} by
choosing some
$q$-deformed reference function (the simplest being the $q$-gaussian
$exp_q[-a|x|^2]$) and
by imposing validity of the $q$-deformed Stoke's theorem
\be
\int d_qV \p^i f(x)=0,~~~~~~~~~~~~~~~~~~f\in\FR.
\ee
$q$-integration, as Riemann one, also satisfies some other important
properties, namely
the reality condition
\be
(\int d_qV~f)^*=\int d_qV~f^*
\ee
for any q ($\in \rn^+$), and
the positivity condition
\be
\int d_qV~f^*f  \ge 0,~~~~~~~~~~~~~~\int d_qV~f^*f =0 \Leftrightarrow f=0;
\ee
finally, it reduces to Riemann integration in the limit $q\rightarrow 1$,
since the
abovementioned algebraic recursion reduce to the classical one in the same
limit.

\sect{Time-independent quantum mechanics on $\rn_q^N$: the q-isotropic
harmonic oscillator}

{}~~~~The formal tools briefly presented
in the preceding section are the ``bricks'' which we use to formulate
time-independent quantum mechanics of one-particle systems on $\rn_q^N$
configuration space.
The idea is to formulate
the time-independent Schroedinger equation as a $q$-differential equation by
introducing a
$q$-deformed hamiltonian $h\in\DFR$ such as
\be
h=-\Delta+ V(x)
\ee
and writing it as an eigenvalue equation
\be
h\psi=E\psi~~~~~~~~~~~~~~~\psi\in\FR.
\ee
A serious problem arises when requiring the hamiltonian to be a hermitean
operator, since
the above laplacian is not real. We solved it in a $nonstandard$ way in Refs.
\cite{fio2,fioth,fioeu}
in two concrete models by postulating that each abstract vector $|u>$ of
the Hilbert space ${\cal H}$
of states can be realized in two different ways $(\psi_u,\bar\psi_u)$
as a $q$-deformed wave-function (with a linear bijective map
$\psi_u\leftrightarrow\bar\psi_u$),
and correspondingly each abstract operator ${\cal B}$ on ${\cal H}$
in two different ways $(b,\bar b)$ as a differential operator. $b$ acts on
$\psi_u$ and
$\bar b$ on $\bar\psi_u$, in other words
${\cal B}|u>:=(b\psi_u,\bar b\bar\psi_u)$;
$\psi_u,b$ (resp. $\bar\psi_u,\bar b$)  will be called the unbarred (resp.
barred)
realization of $|u>,{\cal B}$. Of course this is done in such a way that all
the physics
(eigenvalues of observables, etc.) are the same in either realization.
The crucial point of this approach is the definition of the scalar product of
two vectors
$|v>,|u> \in {\cal H}$ as the sum of the two `` conjugate '' terms:
\be
<u|v>:=\int d_qV\bar{\psi}_u^*\psi_v + \int d_qV \psi^*_u \bar \psi_v .
\ee
Indeed $<~~|~~>$ is manifestly sesquilinear and
\be
<v|u>^*=(\int d_qV\bar{\psi}_v^*\psi_u)^* + (\int d_qV \psi^*_v \bar \psi_u)^*
\stackrel{(1.14)}{=}\int d_qV\psi_u^*\bar{\psi}_v + \int d_qV \bar{\psi}^*_u
\psi_v =<u|v>.
\ee
Its positivity has to be shown for each specific choice of $h$ in formula
(2.1);
 for the proof
in the case of the harmonic oscillator see Ref. \cite{fio2,fioth}. The scalar
product (2.3)
is such that it allows to define a formally hermitean kinetic part $pCp$ of the
abstract hamiltonian $H=pCp+V$ by realizing it as the $pair$ of conjugated
laplacians
$(-q^N\Delta,-q^{-N}\bar{\Delta})$ (here $\bar{\Delta}:=\bp C\bp$).
In fact, it is easy to check that the hermitean conjugate
${\cal B}^{\dagger}$ of an abstract operator ${\cal B}:=(b,\bar b)$ w.r.t. the
scalar product (2.3) is given by the rule
\be
{\cal B}^{\dagger}=(\bar b^*,b^*),
\ee
and consequently $(pCp)^{\dagger}=pCp$, since $\Delta^*=q^{-2N}\bar{\Delta}$.

{}~~~Now we sketch how the previous program can be successfully developed in
the
case of the
$Fun(SO_q(N,\rn))$-isotropic harmonic oscillator. We build $\forall q\in \rn^+$
a sensible quantum mechanical model $\rn^N_q$ ($N\ge 3$)
as the simplest q-deformation of the (time-independent)
classical isotropic harmonic oscillator on ordinary
$\rn^N$; correspondingly, the symmetry group
$SO(N,\rn)$ of rotations of the hamiltonian
is deformed into the quantum group symmetry $SO_q(N,\rn)$. The hamiltonian
 has a lower bounded energy
spectrum and  the scalar product is strictly
positive for any $q\in \rn^+$.
Generalizing the classical algebraic construction, the Hilbert
space of physical states is built applying creation operators to the
(unique) ground state. Observables will be defined as hermitean
operators, as usual. In particular we  construct the observables
hamiltonian, square angular momentum, angular momentum components,
position operators, momentum operators. As in the classical case, the first two
and any angular momentum component will commute with each other; when $N=3,4$
they make up a complete set of commuting observables.

As a first task we have to
fix within the differential algebra $Dif\!f(\rn_q^N)$
a suitable q-analogue $h_{\omega}$ of the classical hamiltonian
$h^{cl}_{\omega}:=-\Delta+\omega^2 xCx$
($x,\p$ being classical coordinates and derivatives)
with characteristic constant $\omega$.
{\it A priori} we don't require it to be necessarily of the form (2.1).
Minimal requirements on $h_{\omega}\in Dif\!f(\rn_q^N)$
are of course that:
\begin{itemize}
\item 1) it should be
a $\FS$-scalar (this is the meaning of the word `` isotropic '' in the
q-deformed setting);

\item 2) it should have a homogenous natural dimension
$d(h_{\omega})=d(\p^2)=2$ and should
reduce to $h^{cl}_{\omega}$ in the limit $q\rightarrow 1$;

\item 3) it should be the ``unbarred'' configuration-space realization of a
hermitean
operator $H_{\omega}$ on some Hilbert space ${\cal H}$ (to be defined);

\item 4) the spectrum of $H_{\omega}$ in ${\cal H}$
should be bounded from below, in order that
$H_{\omega}$ can be considered as the hamiltonian of a sensible
(i.e. stable) quantum mechanical system.

\end{itemize}

Convenience suggests two further requirements.
$h_{\omega}$ being a scalar, it commutes with $U_q(so(N))$ and in particular
with the square angular momentum  $l\cdot l:= l^{ij}l_{ji}$;
this means that when realized as operators,
$h_{\omega}$,$l\cdot l$ can be diagonalized simultaneously.
Requirements 1) - 4) still leave a great freedom in defining
$h_{\omega}$, which can be essentially summarized as its yet undefined
$l\cdot l$-dependence.
By requiring that, as it happens for the classical isotropic harmonic
oscillator,
\begin{itemize}
\item 5) Energy levels~are~$(l\cdot l)$-independent
\end{itemize}
we essentially impose a trivial dependence of $h_{\omega}$ on $l\cdot l$.
Moding out an essential
dilaton-dependence and with a careful choice of the numerical coefficients,
the final unbarred realization of the hamiltonian turns out to be
\be
h_{\omega}:=(-q^N\Delta+\omega^2xCx).
\ee
We introduce the barred hamiltonian
$\bar h_{\omega}:=(-q^{-N}\bar{\Delta}+\omega^2xCx)$ in such a
way that $h_{\omega}^*=\bar h_{\omega}$ and we can build a formally hermitean
abstract hamiltonian
$H_{\omega}$:
\be
H_{\omega}:=(h_{\omega},\bar h_{\omega})~~~~~~~~~~\Rightarrow~~~~~~~~
H_{\omega}^{\dagger}=H_{\omega}.
\ee

Finally, we add the important requirement
\begin{itemize}
\item 6) we would like to algebraically solve the Schroedinger equation
through the introduction of $\FS$-vectors of energy creators
$\vec{A}^+:=(A^{i+})\in Dif\!f(\rn_q^N)$
and annihilators $\vec{A}:=(A^i)\in Dif\!f(\rn_q^N)$,
as in the classical case. More precisely, we require commutations
relations of the type
\be
H_{\omega}A^{i\pm}=A^{i\pm}f^{i\pm}(H_{\omega}),~~~~~~~~~~~f^{i\pm}
(t)\in \cn[t].
\ee
\end{itemize}

In fact, if  relations of the form (2.8) are satisfied, given an eigenvector
$|u>$ of $H_{\omega}$, $H_{\omega}|u>=E|u>$, then $A^{i\pm}|u>$
would be an eigenvector of $H_{\omega}$ with eigenvalue $f^{i\pm}(E)$.

\begin{prop}
A solution of the above problem (2.8) is given by
\be
A^{i\pm}=(a^{i\pm},\bar a^{i\pm}),~~~~~~~~~~~~~~~~\cases{
a^{i\pm}=\La^{-\f 12}[x^i\beta^{\pm}(q,h_{\omega})+\p^i
\gamma^{\pm}(q,h_{\omega})]\cr
\bar a^{i\pm}=\La^{\f 12}[x^i\beta^{\pm}(q^{-1},h_{\omega})+
\p^i\gamma^{\pm}(q^{-1},h_{\omega})] \cr}
\ee
\footnote{In Ref. \cite{fio2} the operators $a^{i\pm}$ were introduced in
block-form as the
collection $\{a^{i\pm}_r\}$ of their projections on the $r$-th eigenspace of
$h_{\omega}$
(see eq. (2.29) below); the here presented more elegant form (where index $r$
is replaced by the
dependence on $h_{\omega}$) was first introduced in Ref. \cite{wata}.}
where $\beta^{\pm},\gamma^{\pm}$ satisfy the condition
\be
\f{\beta^{\pm}(q,h_{\omega})}{\gamma^{\pm}(q,h_{\omega})}=
\f{q^{-N}}{1+q^{2-N}}(h_{\omega}-qf^{\pm}).
\ee
and $f^{i\pm}$ coincide with one of the two solutions $f^{\pm}(h_{\omega})$
of the equation
\be
(qh_{\omega}-f)(q^{-1}h_{\omega}-f)=(1+q^{2-N})^2q^{N-2}\omega^2.
\ee
\end{prop}
It is easy to check that in the limit $q\rightarrow 1$ both $a^{i+}$ and
$\bar a^{i+}$
(respectively $a^{i-}$ and $\bar a^{i-}$) go to the classical creation (resp.
annihilation)
operators of the ordinary isotropic harmonic oscillator.

Now let $\{|u>,E\}$ denote the pair consisting of a eigenvector
$|u>=(\psi_u,\bar\psi_u)$
of $H_{\omega}$ and the corresponding eigenvalue $E$. We can now generate
a ``shower'' of such pairs by means of iterated use of relations (2.8), (2.11):
\be
\matrix
{&     &      &          &             &\{A^{j+}A^{i+}|u>,f^+(f^+(E))\}&...&
\cr
 &  &  &\{A^{i+}|u>,f^+(E)\} & {\nearrow\atop\rightarrow}
                                       &\{A^{j-}A^{i+}|u>,f^-(f^+(E))\}&...&
\cr
 &\{|u>,E\}&{\nearrow \atop \searrow}&  &     &                        &...&
\cr
 &  &  &\{A^{i-}|u>,f^-(E)\} & {\rightarrow\atop\searrow}
                                      &\{A^{j+}A^{i-}|u>,f^+(f^-(E))\}&...& \cr
 &  &   &        &                    &\{A^{j-}A^{i-}|u>,f^-(f^-(E))\}&...&
\cr}
\ee
As in the classical case ($q=1$), starting from an arbitrary $\{|u>,E\}$ we
would generally get
an unbounded (from below) energy spectrum $\{E,f^+(E),f^-(E),...\}$
(in the case $q\neq 1$ it would also be uncountable), which is
non-physical. This is excluded in the case $q=1$ by the condition that $|u>$
is normalizable; this is equivalent to the condition that the energy spectrum
generated starting
from $\{|u>,E\}$ is bounded from below.
It is quite difficult to impose normalizability at this stage in the case
$q\neq 1$,
so we impose directly the second condition.
This implies that one of the generated eigenvector is a ``ground-state'' $|0>$:
\be
A^{i-}|0>=0,~~~~~~~~~~~~~~~H_{\omega}|0>=E_0|0>,
\ee
where
\be
E_0=\omega(q^{{N\over2}-1}+q^{1-{N\over2}})[{N\over 2}]_q, ~~~~~~~~~
[x]_q:=\f{q^x-q^{-x}}{q-q^{-1}}.
\ee
As a direct consequence, for the set of eigenvectors generated by the shower
(2.8) applied to $|0>$
it is true that $\{E=f^+(f^-(E))=f^-(f^+(E))$ i.e. that
$|u>,A^{j-}A^{i+}|u>,A^{j+}A^{i-}|u>$ have the same energy (which was not true
in general),
so that the spectrum is countable (which is necessary for the separability of
the Hilbert space)
and  $A^{i+},A^{i-}$ actually create and destroy an energy excitation
respectively.

Let us call ${\cal H}$ the linear span of all the vectors generated by the
shower (2.8) starting from $|u>=|0>$.
Let us now define
\be
|i_r,i_{r-1},...i_1>:=A^{+i_r}A^{+i_{r-1}}...A^{+i_1}|0>
\ee
\be
{\cal H}_r:=Span_{\cn}\{|i_r,i_{r-1},...i_1>,~~~i_j=-n,...,n\}~~~~~~~~~~~~~~~~
n\ge 0
\ee
\begin{prop}
$r\neq s$ implies ${\cal H}_r\bot {\cal H}_s$ w.r.t. the scalar product (2.4);
 ${\cal H}_r$ has the same dimension
${N+r-1 \choose N-1}$ as in the classical case and is an eigenspace of
$H_{\omega}$ with eigenvalue $E_r$:
\be
H_{\omega}{\cal H}_r=E_r{\cal H}_r~~~~~~~~~~~E_r=
\omega(q^{{N\over2}-1}+q^{1-{N\over2}})[{N\over 2}+r]_q;
\ee
$A^{i\pm}({\cal H}_r)={\cal H}_{r\pm 1}$.
${\cal H}=\bigoplus \limits_{r=1}^{\infty}{\cal H}_r$. Finally, ${\cal H}$
endowed with the scalar
product (2.4) is a pre-Hilbert space, and can be completed into a Hilbert
space in the standard way.
\end{prop}
{\bf Remark} Note the $q\rightarrow q^{-1}$ invariance of the energy levels
$E_r$.  Both for $q>1$ and for
$0<q<1$ the difference $E_{r+1}-E_r$ increases and {\it diverges} with r,
implying that it becomes a macroscopic energy
gap for sufficiently large $r$ (contrary to what happens with some other
q-deformed harmonic
oscillators \cite{bie}).

The definition of $A^{i\pm}$ we gave in equation (2.9) is not complete yet,
since equation (2.10) fixes only the ratio
$\f{\beta^{\pm}(q,h_{\omega})}{\gamma^{\pm}(q,h_{\omega})}$, not $\beta^{\pm}$
itself. One can choose the latter in such a way that annihilators/creators are
hermitean conjugate of each other
\be
(A^{i+})^{\dag}=A^{l-}C_{li}~~~~~~~~~~(A^{i-})^{\dag}=A^{l+}C_{li};
\ee
this condition still leaves some arbitrariness in the definition of
$\beta^{\pm}$, amounting
to the normalization of $A^{i+}A^{j-}$ (in Ref. \cite{fio2} we removed it by
the somewhat
arbitrary condition that the operator $a^{i+}+a^{i-}$ didn't depend on $\p$;
any other choice
would be legitimate).

Commutation relations of the following form hold for the creation/annihilation
operators:
$A^{i+}A^{j-}$ ($B,l^{ij}$ were defined in eqs. $(1.9),(1.10)$; the explicit
rather lenghty form of
the functions $g_A,g_1,g_1',g_S$ is not necessary here):
\be
{\cal P}_{A~hk}^{~~ij}A^{h-}A^{k-}=0={\cal P}_{A~hk}^{~~ij}A^{h+}A^{k+}
\ee
\be
{\cal P}_{A~hk}^{~~ij}A^{h+}A^{k-}=g_A(h_{\omega})l^{ij}~~~~~~~~~~~
{\cal P}_{A~hk}^{~~ij}A^{h-}A^{k+}=-g_A(h_{\omega})l^{ij}
\ee
\be
A^{h+}A^-_h=g_1(h_{\omega})B~~~~~~~~~~~A^{h-}A^+_h=g_1'(h_{\omega})B
\ee
\be
{\cal P}_{S~hk}^{~~ij}A^{h+}A^{k-} - g_S(h_{\omega}){\cal P}_{S~hk}^{~~ij}
A^{h-}A^{k+1}=0
\ee

Due to relation (2.18), $X^j:=A^{j+}+A^{j-}$, $P^j:={\f 1i}(A^{j+}-A^{j-})$
have the same
hermitean conjugation relations of $x^i$, therefore there exist $N$ linearly
independent
combinations of the $X^i$'s (resp. $P^i$'s) which are hermitean operators.
The latter can be called the position and momentum operators respectively,
since in the limit
$q\rightarrow 1$ they become the classical position and momentum operators
respectively.
They don't commute, rather ``$q$-commute'':
\be
{\cal P}_{A~hk}^{~~ij}X^h X^k =0.     ~~~~~~~~~~~~~
{\cal P}_{A~hk}^{~~ij}P^h P^k =0.
\ee

Since $h_{\omega},\bar h_{\omega}$ are scalars, $[l^{ij},h_{\omega}]=0=
[l^{ij},\bar h_{\omega}]$.
One can easily show that angular momentum operators coincide in the barred and
unbarred scheme,
so we can rewrite the above equations in an abstract form
$[l^{ij},H_{\omega}]=0$ where now
by $l^{ij}$ we mean in fact ( with a slight abuse of notation)
the pair $(l^{ij},l^{ij})$; on these
pairs hermitean conjugation coincides with complex conjugation. $l\cdot l$ is a
 scalar itself
and therefore $[l\cdot l,l^{ij}]=0$; moreover, it is hermitean, since it is
real.
Then we can diagonalize simultaneously $H_{\omega}, l\cdot l$ and a real
Cartan subalgebra of $U_q(so(N))$.
This is the q-deformed analogue of diagonalizing angular momentum observables
and a  hamiltonian with central potential in the classical case.
When $N=3,4$ these operators will make up a complete set of commuting
observables.
The eigenvalues of $l\cdot l$ are
\be
l^2_k=[k]_q[k+N-2]_q \f{(q^{2-\f N2}+q^{\f N2 -2})}
{(q+q^{-1})(q^{1-\f N2}+q^{\f N2-1})}~~~~~~~~~~~~~k=0,1,2,...
\ee
One can decompose each ${\cal H}_r$ into the direct sum of eigenspaces of
$l\cdot l$ corresponding to different eigenvalues $l^2_k$. Summing up:
\begin{prop}
${\cal H}_{r,r-2m}$ ($r\ge 0$, $0\le m\le {r\over 2}$)
is an eigenspace of the operators $H_{\omega},l\cdot l$
with eigenvalues $E_r$, $l^2_{r-2m}$ (see (2.17),(2.24)) respectively. Moreover
\be
{\cal H}=\bigoplus\limits_{r=0}^{\infty}\bigoplus \limits_{0\le m \le
{r\over 2}}{\cal H}_{r,r-2m}.
\ee
$\bigoplus$ is to be understood in the sense of
direct sum of mutually $orthogonal$
subspaces w.r.t. the scalar product $<~~|~~>$.
\end{prop}

We give some explicit formulae regarding ``wavefunctions''. The unbarred and
barred solutions
$\psi_0,\bar{\psi}_0\in \FR$ of equation (2.13) are
\be
\psi_0:=e_{q^2}[-{q^{-N}\omega xCx \over 1+q^{2-N}}]~~~~~~~~~~
\bar{\psi}_0:=e_{q^{-2}}[-{q^N\omega xCx \over 1+q^{N-2}}];
\ee
according to our conventions $|0>=(\psi_0,\bar{\psi}_0)$. Introducing the
notation
\be
|i_r,i_{r-1},...i_1>=:(\psi_r^{i_ri_{r-1}...i_1},
\bar {\psi}_r^{i_ri_{r-1}...i_1}),
\ee
one finds
\be
\psi_r^{i_r i_{r-1} ...i_1}:=a^{i_r+}_ra_{r-1}^{i_{r-1}+}...a_1^{i_1+}\psi_0
{}~~~~~~~~~~~~
\bar{\psi}_r^{i_r i_{r-1} ...i_1}:=\bar a^{i_r+}_r\bar a_{r-1}^{i_{r-1}+}...
\bar a_1^{i_1+}\bar{\psi}_0,
\ee
where
\be
a_h^{i+}:=b_h(q)(x^i-{q^{2-h}\over \omega}\p^i)\La^{-\f 12}~~~~~~~~~~~
\bar a_h^{i+}:=b_h(q^{-1})(x^i-{q^{h-2}\over \omega}\bp^i)\La^{\f 12}
{}~~~~~~~~~i=1,2,...,N;
\ee
the pair $(a_h^{i+},\bar a_h^{i+})$ is the explicit form of
$A^{i+}|_{{\cal H}_{h-1}}$,
the creation operator $A^{i+}$ with domain of definition restricted to
${\cal H}_{h-1}$.
Similarly, $A^{i-}|_{{\cal H}_{h-1}}=:(a_h^{i-},\bar a_h^{i-})$
where
\be
a_h^i:=q^{2-2h-N}b_h(q)(x^i+{q^{h+N}\over \omega}\p^i)\La^{-\f 12}~~~~~~~~~
\bar a_h^i:=q^{-2+2h+N}b_h(q^{-1})(x^i+{q^{-h-N}\over \omega}\p^i)\La^{\f 12};
\ee
$b_h,d_h$ are given in Ref. \cite{fio2,fioth}. One can easily verify from
equations (2.28) that
\be
\cases{\psi_r^{i_r i_{r-1} ...i_1}=P_r(x)e_{q^2}[-{q^{-N-r}\omega xCx \over
(1+q^{2-N})}]\cr
\bar {\psi}_r^{i_ri_{r-1}...i_1}=P_r'(x)e_{q^{-2}}[-{q^{r+N}\omega xCx \over
(1+q^{N-2})}] \cr}
\ee
where $P_r(x),P'_r(x)$ are two polynomials of degree $r$ in $x$ containing
only terms of degree
$r, r-2, r-4,...$, and $e_q[Z]$ denotes the q-exponential
\be
e_q[Z] :=\sum\limits_{n=0}^{\infty}{Z^n \over (n)_q!},
{}~~~~~~~~~(n)_q:=\f{q^n-1}{q-1}.
\ee
$\psi_r^{i_ri_{r-1}...i_1}, \bar {\psi}_r^{i_ri_{r-1}...i_1}$ can be called
the unbarred,barred
q-deformed Hermite functions, since they both reduce to the classical Hermite
functions in the limit $q\rightarrow 1$.

Finally, the wavefunctions realizing the subspace
${\cal H}_{r,r-2m}$ ($r\ge 0$, $0\le m\le {r\over 2}$) are combinations of the
$\psi_r^{i_ri_{r-1}...i_1}, \bar {\psi}_r^{i_ri_{r-1}...i_1}$ with fixed $r$.
In the unbarred case
they are given by
\be
[({\cal P}_1\ot'...\ot' {\cal P}_1\ot' {\cal P}_{r-2m,S})\psi_r]
^{l_1...l_r}\propto
P^{l_{2m+1}...l_r}_S(x) p_{r,m}(xCx)
e_{q^2}[-{\omega q^{-r-N}xCx\over \mu}],
\ee
where $p_{r,m}(z)$, $ 0\le m\le {r\over 2}$ are polynomials in one variable $z$
 and
\be
P^{l_{2m+1}...l_n}_S(x):={\cal P}_{r-2m,S,j_{2m+1}...j_r}
^{~~~~~~~l_{2m+1}...l_r} x^{j_{2m+1}}...x^{j_r}
\ee
are the ``q-deformed spheric homogenous polynomials'' of degree $r$ obtained
through application
of the q-symmetric-modulo-trace projector ${\cal P}_{r-2m,S}$. The latter is
defined by
\be
{\cal P}_{k,S}{\cal P}_{Ai,(i+1)}=0={\cal P}_{k,S}{\cal P}_{1i,(i+1)},
{}~~~~~~~~~
{\cal P}_{Ai,(i+1)}{\cal P}_{k,S}=0={\cal P}_{1i,(i+1)}{\cal P}_{k,S},
{}~~~~~~~~~~~~({\cal P}_{k,S})^2={\cal P}_{k,S},
\ee
$1\le i \le k-1,$
where ${\cal P}_{Ai,(i+1)}=(\otimes {\bf 1})^{i-1}\otimes{\cal P}_A\otimes
(\otimes {\bf 1})^{n-i-1}$, etc. The formulae in the barred case are similar.
For further details on wavefunctions see also U. $\&$ S. Watamura in these
proceedings.

\section*{Acknowledgments}
I thank V.K. Dobrev and H.-D. Doebner for inviting me to the symposium and J.
wess for his kind
hospitality at his Institute.

\end{document}